%% file: IRO24.tex
\documentclass[12pt,oneside,openany,article]{memoir}
\usepackage{mempatch}

\nouppercaseheads 
\usepackage[amsthm,thmmarks,hyperref]{ntheorem}
\usepackage[noTeX]{mmap}
\usepackage{etex}

\usepackage[english]{babel}
\usepackage[leqno]{mathtools}

\usepackage{mathrsfs}
\usepackage{upgreek}
\usepackage[compress,square,comma,numbers]{natbib}
\usepackage{hypernat} %required if hyperref is used as well

\usepackage{sfmath}

\usepackage{microtype}

\usepackage{subfig}
\usepackage{wrapfig}
\usepackage{array}

\usepackage{graphicx,color}
\definecolor{gray}{gray}{0}
\pagecolor{white}

\usepackage{hyperxmp}
\usepackage{xr-hyper}
\usepackage{nameref}
\usepackage[pdftex,bookmarks,pdfnewwindow,plainpages=false,unicode]{hyperref}

\usepackage{bookmark}

\usepackage{enumitem}

\usepackage{amssymb} % for \Subset, \lozenge, etc.

\hypersetup{
colorlinks=true,
linkcolor=black,
citecolor=black,
urlcolor=blue,
pdfauthor={Victor Ivrii},
pdftitle={Operators with Periodic Hamiltonian Flows in Domains with the Boundary},
pdfsubject={Sharp Spectral Asymptotics},
pdfkeywords={Microlocal Analysis,  Sharp Spectral Asymptotics, Periodic flows},
bookmarksdepth={4}
}

\numberwithin{equation}{chapter}

\input{defines_IRO24}

\begin{document}
\title{Operators with Periodic Hamiltonian Flows in Domains with the Boundary}
\author{Victor Ivrii}

\maketitle
{\abstract%
We consider operators in the domains with the boundaries and derive sharp spectral asymptotics (containing non-Weyl correction) in the case  when Hamiltonian flow is periodic. 

Even if operator is scalar but not second order (or even second-order but there is an inner boundary with both refraction and reflection present) Hamiltonian flow is branching on the boundary and the notion of periodicity becomes more complicated.

\endabstract}

\setcounter{chapter}{-1}
\chapter{Introduction}
\label{book_new-sect-5-0}

In section \ref{book_new-sect-6-2} of \cite{futurebook}\footnote{\label{foot-15-0} This article is a rather small part of the huge project to write a book and is just section~8.3 consisting entirely of newly researched results. section~\ref{book_new-sect-6-2} roughly corresponds to Section~4.7 of its predecessor V.~Ivrii~\cite{Ivr1}. External references by default are to  \cite{futurebook}.}

we considered the case when there is no boundary and Hamiltonian flow is periodic and derive very sharp spectral asymptotics. Now we want to achieve similar results when there is a boundary. We do not expect asymptotics to be that sharp, but better than $O(h^{1-d})$ and have non-Weyl correction term.  However the presence of the boundary, more precisely, branching of the rays on the boundary brings the new possibilities. 

Even if operator is scalar but not second order (or even second-order but there is an inner boundary with both refraction and reflection present) Hamiltonian flow is branching on the boundary and the notion of periodicity becomes more complicated.

For simplicity we consider only Schr\"odinger operators and derive sharp spectral asymptotics (containing non-Weyl correction) in the case  when Hamiltonian flow is periodic. 

\chapter{Discussion and plan}
\label{sect-8-3-1}

We start our analysis in subsection \ref{sect-8-3-2} from the case when there is no branching and the Hamiltonian flow (with reflections) is periodic.  Then if period is $T_0$, $\Psi_{T_0}=I$ we have the same equality 
\begin{equation}
e^{ih^{-1}T_0A}Q= e^{i T_0 L}Q
\label{8-3-1}
\end{equation}
as before where $\supp q$ is disjoint from the boundary and only transversal to the boundary trajectories originate at $\supp q$ and $L$ is $h$-pseudo-differential operator. Then 
\begin{gather}
F_{t\to h^{-1}\tau} \bar{\chi} (\frac{t}{T_0}-n) 
\Gamma \bigl(Q_{1x}u\,^t\!Q_{2y}\bigr)
=(2\pi h)^{1-d}J(n)+O(h^{1-d+\delta})
\label{8-3-2}\\
\shortintertext{with}
J(n)=\int _{\Sigma_\tau} e^{in \ell}q_1q_2 d\upmu_\tau 
\label{8-3-3}
\end{gather}
as $\bar{\chi}$ is supported in $(-\frac{2}{3},\frac{2}{3})$ and equals $1$ in 
$(-\frac{1}{3},\frac{1}{3})$ and $n\in \bZ$, $|n|\le h^{-\delta}$. Here as before $\Sigma_\tau=\{a(x,\xi)=\tau\}$ and $\upmu_\tau =dx d\xi:da$ is the measure there, $\ell$ is the principal symbol of $L$.

If $\ell$ is ``variable enough'' on $\Sigma_\tau$ so the right-hand expression is decaying as $n\to \infty$ we can improve remainder estimate $O(h^{1-d})$ (to what degree depends on the rate of the decay). In section \ref{book_new-sect-6-2} essentially $\ell$ was defined by integral of the subprincipal symbol along periodic trajectories. Now however $\ell $ can pick up extra terms at the moment of reflection: for example, for Schr\"odinger operator  increment of $\ell$ is $0$ and $\pi $ for Neumann and Dirichlet boundary condition respectively; however more general boundary condition brings increment which depends on the point $(x',\xi')\in T^*Y$ of reflection.

Then in section \ref{sect-8-3-3} we consider a more complicated case when two operators are intertwined through boundary condition and all Hamiltonian trajectories of one of them are periodic but of another is not and, moreover, if the generic Hamiltonian billiard is periodic, it is in fact Hamiltonian billiard of the first operator. Then (\ref{8-3-1}) fails but (\ref{8-3-2})--(\ref{8-3-3}) is replaced by (\ref{8-3-2}), (\ref{8-3-4})
\begin{equation}
J(n) =\int _{\Sigma_{1,\tau}} e^{in \ell- |n| \ell '}q_1q_2 d\upmu_{1.\tau}
\label{8-3-4}
\end{equation}
where $\Sigma_{1,\tau}$ and $\upmu_{1,\tau}$ correspond to the first operator and $\ell'\ge 0$ depends on the portion of energy which was whisked away along trajectory of the second operator at reflection point. Then if 
\begin{equation}
\upmu_{1,\tau}\bigl(\{(x,\xi)\in \Sigma_{1,\tau}, \ell'=0, 
\nabla \ell=0\}\bigr)=0
\label{8-3-5}
\end{equation}
we conclude that $J(n)=o(1)$ as $n\to \infty$ and we will be able to improve remainder estimate $O(h^{1-d})$ to $o(h^{1-d})$ and again there will be a non-Weyl correction term but it will depend also on $\ell$ and $\ell'$.

Finally, in subsection  \ref{sect-8-3-3} we consider even more complicated when two operators are intertwined through boundary condition and Hamiltonian trajectories of both of them are periodic and moreover after a while all trajectories issued from some point at $T^*Y$ assemble there. Then we can construct a matrix symbol and its eigenvalues  play a role of symbol $\ell$.

\chapter{Simple Hamiltonian flow}
\label{sect-8-3-2}

\section{Inner asymptotics}
\label{sect-8-3-2-1}
Let $\cU$ be an open connected subset in $T^*X$ disjoint from $\partial X$. Consider Hamiltonian billiards issued from $\cU$, let $\Phi_t$ denote generalized Hamiltonian flow. Assume that
\begin{gather}
|\nabla_{x,\xi}a(x,\xi)|\ge \epsilon \qquad \forall t\ \ 
\forall (x,\xi)\in \Phi_t(\cU), \label{8-3-6}\\[3pt]
|\phi(x)| + |\{a,\phi\}(x,\xi)|\ge \epsilon\qquad \forall t\ \ 
\forall (x,\xi)\in \Phi_t(\cU),\label{8-3-7}\\[3pt]
\Phi_t (x,\xi)=(x,\xi)\qquad \text{with\ \ }t=t(x,\xi)>0\ \ \forall (x,\xi)\in \cU,\label{8-3-8}
\end{gather}
where $\phi\in \sC^\infty (\bar{X})$, $\phi >0$ in $X$ and $\phi =0$ on $\partial X$ and 
\begin{claim}\label{8-3-9}
Along $\Phi_t(\cU)$ reflections are simple without branching, i.e. 
$\iota ^{-1}\iota(y,\eta)\cap\{a(y,\eta)=\tau\}$
consists of two (disjoint) points $(y,\eta^\pm)$ such that $\pm \{a,\phi\}>0$
as $(y,\eta)\in \Phi_t(\cU)$, $y\in \partial X$ and $\tau=a(y,\eta)$.
\end{claim}
Let us recall that $\iota:T^*X|_{\partial X}\to T^*\partial X$ is a natural map.

Then exactly as in the case without boundary
\begin{claim}\label{8-3-10}
One can select $t(x,\xi)= T(a(x,\xi))$ in $\cU$ with $T\in \sC^\infty$  (so period depends on energy level only\footnote{\label{foot-8-6} But there could be subperiodic trajectories with periods $n^{-1}T(a(x,\xi))$, $n=2,3,\dots$. The structure of subperiodic trajectories described in subsection~\ref{book_new-sect-6-2-5}.} 
\end{claim}
and
\begin{claim}\label{8-3-11}
As $(x,\xi)\in \cU$ and $0<t< T(x,\xi)$ \ $\Phi_t(x,\xi)$ meats $\partial X$ exactly $(N-1)$ times with $N=\const$; so as $N\ge 1$ billiard consists of $N$ segments with the ends on the boundary.
\end{claim}

We can consider instead of $A$ with domain 
\begin{equation}
\fD(A)=\{u\in \sH^m(X), \eth Bu=0\}
\label{8-3-12}
\end{equation}
operator $f(A)$ with  
\begin{equation}
f(\tau)=\int T^{-1}(\tau)\,d\tau
\label{8-3-13}
\end{equation}
and introduce propagator  $e^{ih^{-1}tf(A)}$
so that the Hamiltonian flow $\Psi_t$ of $f(A)$ is periodic with period $1$.

However now structure of $f(A)$ near $\partial X$ is a bit murky; but it is not really serious obstacle since we can consider problem
\begin{gather}
\bigl( \check{f}(hD_t) - A\bigr) e^{ih^{-1}tf(A)}=0,\label{8-3-14}\\[3pt]
\eth B  e^{ih^{-1}tf(A)}=0 \label{8-3-15}
\end{gather}
where $\check{f}$ is an inverse function. Therefore its Schwartz kernel $U(x,y,t)$ satisfies problem (\ref{8-3-14})--(\ref{8-3-15}) with respect to $x$ and transposed problem with respect to $y$.

Note that due to assumptions (\ref{8-3-6})--(\ref{8-3-7})  $e^{ih^{-1}tf(A)}Q$ is a Fourier integral operator with the symplectomorphism $\Psi_t$ provided $Q$ is $h$-pseudo-differential operator with the symbol supported in $\cU$. Further, due to  (\ref{8-3-8}), (\ref{8-3-13}) $\Psi_1=I$ and therefore $e^{ih^{-1}f(A)}Q$ is $h$-pseudo-differential operator. So, we arrive to

\begin{proposition}\label{prop-8-3-1}
Let $A$ be a Schr\"odinger operator with the scalar principal symbol $a(x,\xi)$ satisfying \textup{(\ref{8-3-6})}--\textup{(\ref{8-3-8})} and \textup{(\ref{8-3-13})}. Then
\begin{equation}
e^{ih^{-1}f(A)-i\kappa }\equiv e^{ih^{-1}\varepsilon L}\qquad\text{in\ \ }\cU
\label{8-3-16}
\end{equation}
with $\varepsilon =h$ and $h$-pseudo-differential operator $L$; here $\kappa$ is Maslov's constant and the principal symbol $\ell$ of $L$ is defined by
\begin{equation}
e^{i\ell }= \prod_{1\le k\le N} e^{i\ell _k}e^{i\ell' _k}
\label{8-3-17}
\end{equation}
where $\ell_k$ corresponds to $k$-th segment (see theorem \ref{book_new-thm-6-2-3}) and $\ell'_k$ corresponds to $k$-th reflection. 
\end{proposition}

\begin{example}\label{ex-8-3-2}
Consider Schr\"odinger operator. Assume that at the reflection point $a(x,\xi)=\xi_1^2+a(x',\xi')$, $b=\xi_1+i\beta(\xi')$. Then
\begin{equation}
e^{i\ell_k '}= 
- \bigl((\tau -a')^{\frac{1}{2}}-i\beta \bigr)^{-1}
\bigl(  (\tau -a')^{\frac{1}{2}}+\i\beta\bigr) 
\label{8-3-18}
\end{equation}
with the right-hand expression calculated in the corresponding reflection point.

Taking $\beta=0$ or $\beta=\infty$ (formally) we get $\ell'_k=0$ and $\ell'_k=\pi$  for Neumann and Dirichlet boundary conditions respectively.
\end{example}

Then we can use all the local results of section \ref{book_new-sect-6-2}, assuming that (\ref{8-3-16}) holds with $h^l\le \varepsilon \le h$ and derive local spectral asymptotics with the same precision (up to $O(h^{2-d})$) as there.

\begin{corollary}\label{cor-8-3-3} Let conditions of proposition \ref{prop-8-3-1} be fulfilled. Then asymptotics with the remainder $o(h^{1-d})$ holds provided 
\begin{equation}
\upmu_\tau  \bigl(\{(x,\xi)\in \Sigma_\tau ,\nabla _{\Sigma_\tau}\ell (x,\xi)=0\}\bigr)=0
\label{8-3-19}
\end{equation}
where $\nabla_{\Sigma_\tau}$ means a gradient along $\Sigma_\tau$;

\medskip\noindent
(ii) Asymptotics with the remainder $O(h^{1-d+\delta})$ with small exponent $\delta>0$ holds provided 
\begin{equation}
\upmu_\tau  \bigl(\{(x,\xi)\in \Sigma_\tau , |\nabla '\ell (x,\xi)| \le \varepsilon \}\bigr)=o(\varepsilon^{\delta'})
\label{8-3-20}
\end{equation}
with the small exponent $\delta'>0$;

\medskip\noindent
(iii) Furthermore asymptotics with the remainder $o(h^{1-d})$ holds as \textup{(\ref{8-3-6})}--\textup{(\ref{8-3-7})} are replaced by their non-uniform versions 
\begin{equation}
|\phi(x)| + |\{a,\phi\}(x,\xi)|>0\qquad \forall t\ \ 
\forall (x,\xi)\in \Phi_t(\cU\cap \Sigma_\tau\setminus \Lambda_\tau),\tag*{$\textup{(\ref*{8-3-7})}'$}\label{8-3-7-'}
\end{equation}
where $\upmu _\tau (\Lambda_\tau)=0$.
\end{corollary}

We leave to the reader to formulate statement similar to (iii) but with the remainder estimate $O(h^{1-d+\delta})$.

Usually our purpose is the asymptotics with $Q_{1x}=\psi_1(x)$, $Q_{2y}=\psi_2(y)$ where $\psi_j$ are smooth functions but not vanishing near the $\partial X$. 

\begin{corollary}\label{cor-8-3-4}
Asymptotics with the remainder estimate $o(h^{1-d})$ holds provided conditions \textup{(\ref{8-3-6})}, \ref{8-3-7-'}, \textup{(\ref{8-3-8})} and \textup{(\ref{8-3-19})} are fulfilled in $\cU=T^*X$, $Q_j=\psi_j\in \sC^\infty (\bar{X})$ compactly supported.
\end{corollary}

\begin{proof}
Due to corollary \ref{cor-8-3-3}(iii) contribution of the zone $\{x_1>\varepsilon\}$ to the remainder is $o_\varepsilon (h^{1-d})$ while contribution of zone $\{x_1\le \varepsilon\}$ to the remainder does not exceed $O(\varepsilon h^{1-d})$ as $\varepsilon \ge h$.
\end{proof}

These arguments could be improved to $\varepsilon = h^\delta$ and the remainder estimate could be improved to $O(h^{1-d+\delta})$. We leave the precise statements and arguments to the reader. Our goal is to derive more sharp remainder estimate.

\section{Asymptotics near boundary}
\label{sect-8-3-2-2}

Let us consider $X=\{x_1>0\}$ and $\cU$ open subset in $T^*\bR^{d-1}_{x'} \times [0,\varepsilon]_{x_1}\times\bR_\tau $ where $\varepsilon>0$ is a very small constant. We are interested in 
\begin{equation}
F_{t\to h^{-1}\tau}\bar{\chi}_T(t) \Tr \bigl(e^{ih^{-1}tf(A)}Q\bigr)
\label{8-3-21}
\end{equation}
with $Q=Q (x,hD'_x,hD_t)$ with symbol supported in $\cU$; let us rewrite it as 
\begin{equation*}
F_{t\to h^{-1}\tau}\bar{\chi}_T(t) \Tr (e^{ih^{-1}(t-\bar{t})f(A)}Q'e^{ih^{-1}\bar{t}f(A)})=
F_{t\to h^{-1}\tau}\bar{\chi}_T(t) \Tr \bigl(e^{ih^{-1}tf(A)}Q'\bigr)
\end{equation*}
with $Q'=e^{ih^{-1}\bar{t}f(A)}Qe^{-ih^{-1}\bar{t}f(A)}$.

Note that if (\ref{8-3-6}), (\ref{8-3-7}) are fulfilled, $\varepsilon >0$  and $\bar{t}\ge C\varepsilon$ are small enough then  $Q'=Q'(x,hD_x,hD_t)$ is $h$-pseudo-differential operator with the symbol supported in 
$\{\epsilon' \bar{t}\le x_1\le C'\bar{t}\}$ and we can rewrite 
\begin{equation*}
Q'=Q''(x,hD_x)+ Q'''(x,hD_x,hD_t)\bigl(hD_t-f(A)\bigr);
\end{equation*}
therefore we can rewrite (\ref{8-3-21}) in the same form but with $Q$ replaced by $Q''$. 

Now we can apply all the arguments of the previous subsection and of \ref{book_new-sect-6-2} and derive asymptotics with the remainder estimate as sharp as  $O(h^{2-d})$ provided conditions (\ref{8-3-6})--(\ref{8-3-8}) are fulfilled in the full measure as well as corresponding conditions of section~\ref{book_new-sect-6-2} to $\ell$. Here as usual we lift the points of 
$(y',\eta',\tau)\in T^*\partial X$ to $(y',\eta^\pm)\in \Sigma_\tau$.

So, we need to analyze only near tangent zone $\{x_1\le \varepsilon, |\{a,x_1\}|\le \varepsilon\}$. 

Unfortunately the theory of the manifolds with the boundaries and periodic Hamiltonian (or even geodesic) flow with reflections is completely undeveloped.

Let us assume that (\ref{8-3-6}), \ref{8-3-7-'} and (\ref{8-3-8}) are fulfilled for all $(x,\xi)\in T^*X: \tau_1\le a(x,\xi)\le \tau_2$ with $\tau_1<\tau_2$.
Then trajectories tangent  to $\partial X$ cannot penetrate into $X$ and therefore $x_1=0, \tau_1\le a'(x,\xi')\le \tau_2$ is incompatible with  $a'_{x_1}(x,\xi)< 0$. 

Further, $x_1=0, \tau_1\le a'(x,\xi')\le \tau_2$ is incompatible with $a'_{x_1}(x,\xi)>0$ as well; really otherwise almost tangent to boundary billiards make very small jumps which contradicts to periodicity after $N$ reflections with fixed $N$.

Therefore $\partial X$ is \emph{bicharacteristically flat\/} \index{bicharacteristically flat}
\begin{equation}
x_1=0, \tau_1\le a'(x,\xi')\le \tau_2 \implies a'_{x_1}(x,\xi)= 0.
\label{8-3-22}
\end{equation}
Further, (\ref{8-3-8}) implies that 
\begin{claim}\label{8-3-23}
All solutions of equation
\begin{equation}
z'' +b (\Psi_t(x',\xi'))z=0\qquad\text{with\ \ }b(x',\xi)=\frac{1}{2} a'_{x_1x_1}(0,x',\xi')
\label{8-3-24}
\end{equation}
are $T_0$-periodic for any $(x',\xi')\in T^*\partial X$, $\tau_1\le a'(x,\xi')\le \tau_2$.
\end{claim}
Therefore (\ref{8-3-8}) implies that
\begin{equation}
\rho(x,\xi) \circ \Psi_t \asymp \rho\quad \text{with \ \ }\rho = x_1 +|\{a,x_1\}|\quad\text{as \ \ }
\tau_1\le a(x,\xi)\le \tau_2.
\label{8-3-25}
\end{equation}

Consider zone $\cU_\varepsilon =\{(x,\xi): \rho (x,\xi)\asymp \varepsilon\}$. Blowing up $(x_1,\xi_1)\to \varepsilon^{-1}(x_1,\xi_1)$ we can prove easily that after this
\begin{gather}
\Psi _t (x,\xi)= \bigl(\Psi'_{t,\varepsilon}(x ,\xi) , \Psi''_{(x',\xi',t),\varepsilon}(x_1,\xi_1)\bigr),\label{8-3-26}\\
\shortintertext{with}
\Psi'_{t,\varepsilon}(x,\xi)\sim  \Psi'_{t,0}(x',\xi')+\sum_{n\ge 2}
\varepsilon^n  \Psi'_{t,(n)}(x,\xi),\label{8-3-27}\\
\Psi''_{(x',\xi',t),\varepsilon}(x_1,\xi_1)\sim \sum_{n\ge 1}
\varepsilon^n  \Psi''_{(x',\xi',t),(n)}(x_1,\xi_1)\label{8-3-28}
\end{gather}
where $\Psi'_{t,0}$ is the Hamiltonian flow on $T^*\partial X$ and 
$\Psi''_{t,(1)}$ is the linearized with respect to $(x_1,\xi_1)$ billiard flow (described by (\ref{8-3-24})). One can prove easily that $\Psi'_{t,(n)}$ and
$\Psi''_{(x',\xi',t),(n)}$ are uniformly smooth.

Then in blown-up coordinates (\ref{8-3-16}) holds provided $\varepsilon \ge h^{\frac{1}{2}-\delta}$ because $\hbar=\varepsilon^{-2}h$ and then we can apply all the previous arguments and derive asymptotics; the remainder estimate is as sharp as $O(h^{2-d})$ provided (in not blown-up coordinates)
\begin{equation}
|\nabla _{(x,\xi)} \ell (x,\xi)|\ge \epsilon_0\qquad\text{as \ \ } 
x_1\asymp \rho (x,\xi)\le \epsilon_1;
\label{8-3-29}
\end{equation}
$x_|\asymp \rho (x,\xi)$ means exactly that $x_1\ge \epsilon |\xi_1|$.

Therefore 

\begin{proposition}\label{prop-8-3-5}
Let $A$ be a Schr\"odinger operator with the  principal symbol $a(x,\xi)$ satisfying \textup{(\ref{8-3-6})}--\textup{(\ref{8-3-8})} and \textup{(\ref{8-3-22})}--\textup{(\ref{8-3-23})}. Let the subprincipal symbol and the boundary condition be such that \textup{(\ref{8-3-29})} is fulfilled. All these conditions are supposed to be fulfilled in $\{(x,\xi):\rho (x,\xi)\le \epsilon_0\}$.

Then the contribution of zone $\{h^{\frac{1}{2}-\delta}\le \rho (x,\xi)\le \epsilon_2\}$ to the remainder  is $O(h^{2-d})$.
\end{proposition}

Since the contribution of zone  $\{\rho (x,\xi)\le \varepsilon\}$ to the remainder does not exceed $C\varepsilon^2 h^{1-d}$ and we can take $\varepsilon=h^{\frac{1}{2}-\delta}$ we conclude that the contribution of zone $\{\rho (x,\xi)\le \epsilon_2\}$ to the remainder is $O(h^{2-d-\delta})$. One can also estimate this way the contribution of subperiodic trajectories. Thus we arrive (details are left to the reader) to

\begin{theorem}\label{thm-8-3-6}
Let $A$ be a Schr\"odinger operator with the  principal symbol $a(x,\xi)$ satisfying \textup{(\ref{8-3-6})}--\textup{(\ref{8-3-8})} and \textup{(\ref{8-3-22})}--\textup{(\ref{8-3-23})}. Let the subprincipal symbol and the boundary condition be such that \textup{(\ref{8-3-29})} is fulfilled. 
All these conditions are supposed to be fulfilled in 
$\{(x,\xi):\rho (x,\xi)\le \epsilon_0\}$.

Then contribution of zone $\{\rho (x,\xi)\le \epsilon_2\}$ to the remainder  is $O(h^{2-d-\delta})$ while its contribution to the principal part of the asymptotics is given by the standard two term Weyl formula plus a correction term constructed according to section~\ref{book_new-sect-5-2} for symbol $\ell$.
\end{theorem}

\begin{problem}\label{prob-8-3-7}
Recover $O(h^{2-d})$ remainder estimate. This should not be extremely hard especially if there are no subperiodic trajectories but definitely worth of publishing.
\end{problem}

\begin{example}\label{ex-8-3-8}
Consider Laplacian $h^2\Delta$ on the standard hemisphere 
$\bS^d_+\Def \{|x|=1, x_1>0\}$ in $\bR^{d+1}$ or harmonic oscillator  
$h^2\Delta +|x|^2$ on the standard half-space $\bR^d_+=\{x_1>0\}$.

\medskip\noindent
(i) Consider boundary operator of example \ref{ex-8-3-2}. Then  $N=2$ and $\ell (z,\tau)= \ell'(z,\tau)+\ell'(\varrho(z),\tau)$ where $z= (x',\xi')$, $\varrho(z)$ is the \emph{antipodal point\/} \index{point!antipodal} and $\ell'$ is defined by (\ref{8-3-18}). Then one can express conditions to $\ell$ as conditions to $\beta$; we leave exact statements to the reader;  

\medskip\noindent
(ii) Consider Dirichlet or Neumann boundary conditions. Then as perturbation is $h x_1$  one can calculate easily that $\ell = \kappa \rho +O(\rho^2)$ where $\rho$ is the incidence angle of trajectory and $\kappa>0$. Then condition (\ref{8-3-29}) is fulfilled.
\end{example}

\section{Discussion}
\label{sect-8-3-2-3}

As we mentioned almost nothing is known about manifolds with all billiards closed; we are discussing geodesic billiards but one can ask the same questions about Hamiltonian billiards or Hamiltonian billiards associated with the Schr\"odinger operator; in the two last cases phase space $S^*X$ is replaced by the portion of the phase space $\{\tau_1\le a(x,\xi)\le \tau_2\}$. We assume that the phase space is connected:

\begin{problem}\label{prob-8-3-9}
(i) So far our only examples are a standard hemisphere $\bS^d_+$ and a standard half-space $\bR^d_+$ (see figure below); are any other there really different examples? 
\begin{figure}[h!]\centering
\subfloat[Hemisphere, projected to equatorial plane]{\includegraphics[width=.5\textwidth]{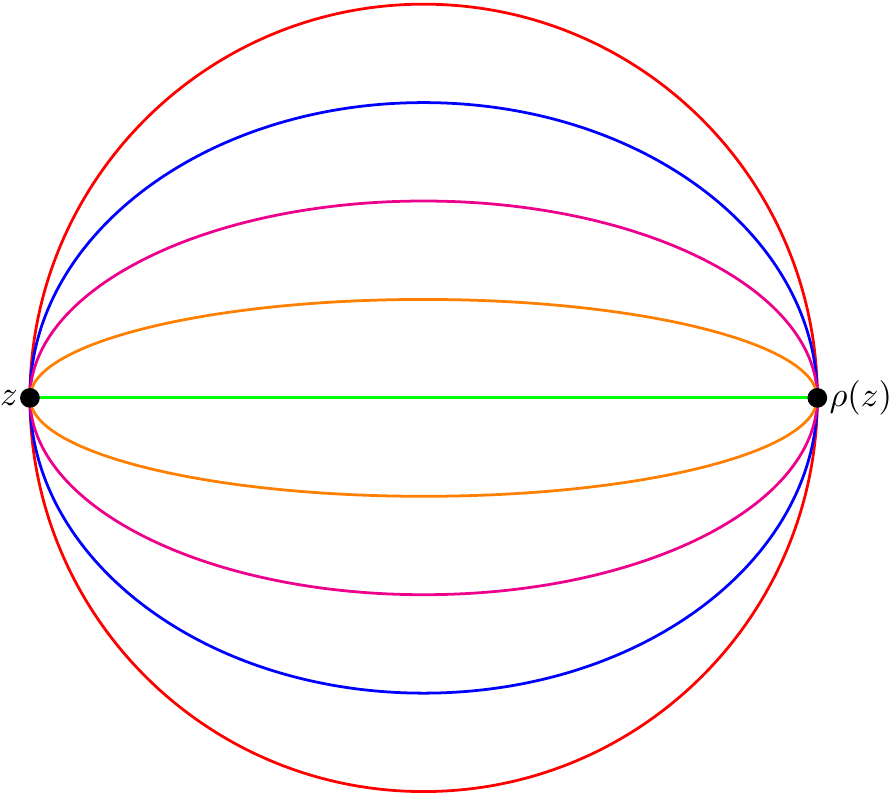}}
\qquad
\subfloat[half-plane]{\includegraphics[width=.25\textwidth]{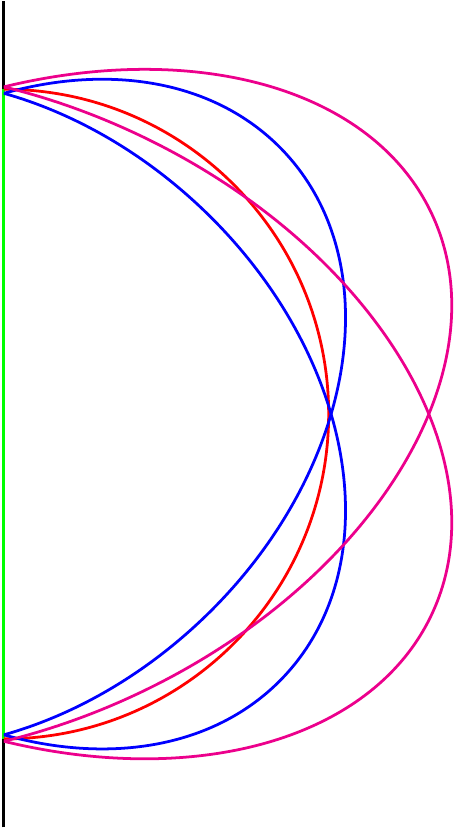}} \caption{\label{fig-billiards} Different billiards (actually, their $x$-projections) are shown by different colors, they reflect from the boundary at two mutually antipodal points.}
\end{figure}

\noindent
(ii) In particular, are there manifolds with $N\ne 2$ (i.e. with $N=1,3,4,\dots$)? We can consider $N$ copies $X_1,\dots, X_N$ of the standard quarter-sphere 
$\bS^2\cap\{0\le\theta\le \frac{\pi}{2}\}\cap\{0\le \phi \le \pi\}$, then glue 
$\{X_k, \phi=\pi\}$ with $\{X_{k+1}, \phi=0\}$ for $k=1,\dots, N$, where $X_{N+1}\Def X_N$. Then for resulting manifold we have $N$ reflections but for $N\ne 2$ there will be singularity at North Pole.

\medskip\noindent
(iii) In particular, are there manifolds with $\partial X$ which is not connected?

\medskip\noindent
(iv) In particular, are there manifolds with subperiodic billiards and (or) with subperiodic boundary trajectories?
\end{problem}

\chapter{Branching Hamiltonian flow with ``scattering''}
\label{sect-8-3-3}

Let us consider two manifolds $X_1$ and $X_2$ glued along the connected component $Y$ of the boundary. We consider Schr\"odinger operators $A_j$ on $X_j$ and these operators are intertwined through the boundary conditions. We assume that Hamiltonian flows have simple (no branching) reflections on each of them, so branching comes from ``reflection-refraction''. Further we are interested in the case when Hamiltonian flow on $X_1$ is periodic while on $X_2$ majority of trajectories are not periodic.

\section{Analysis in $X_2$}
\label{sect-8-3-3-1}

Since we do not know any other examples of manifolds with periodic Hamiltonian flows with reflections but hemisphere or half-space and $N=2$ then we assume that 
\begin{claim}\label{8-3-30}
For $z\in T^*Y$ $\Phi_{1,\tau}(z)=\varrho(z)$~\footnote{\label{foot-8-7} We define $\Phi_{j,\tau}$ in the following way: we lift $z\in T^*Y$ to $z\in T^X_j|_Y\cap \Sigma_{j,\tau}$ with $\{a_j,x_1\}(z)>0$, launch trajectory forward until it hits the boundary at $(y',\eta)\in T^*X_j|_Y$ and then project to $T^*Y$.},  $\varrho$ is an antipodal map on $T^*Y$. 
\end{claim}

We also assume that $\Phi_{1,\tau}$ and $\Phi_{2,\tau}$ commute
\begin{equation}
\Phi_{1,\tau}\circ \Phi_{2,\tau}=\Phi_{2,\tau}\circ \Phi_{1,\tau}.
\label{8-3-31}
\end{equation}
Then  for points of $\Sigma_{2,\tau}$ periodicity properties of the branching Hamiltonian flow with reflections $\Psi_t$ coincide with those of $\Psi_{2,t}$.

Now we need to impose condition to the boundary operator except that $\{A,B\}$ is self-adjoint. Namely, consider a point $z=(x',\xi')\in T^*Y$, and consider an auxiliary problem 
\begin{gather}
a_j(x',D_1,\xi')u_j=0  \qquad j=1,2.\label{8-3-32}\\[2pt]
\eth b(x',D_1,\xi')u=0\label{8-3-33}
\end{gather}
(where $u=(u_1,u_2)$) and consider solutions which are combinations of  exponents $e^{i\lambda_j x_1}$ where for each $j$ \underline{either} 
$\Im \lambda_j>0$ \underline{or} $\Im \lambda_j>0$ and $\{a_j,x_1\}(x',\lambda_j,\xi')>0$. Denote by $\Lambda$ the set of points for which such non-trivial solutions exist. We assume that
\begin{equation}
\mes_{T^*Y}\Lambda=0.
\label{8-3-34}
\end{equation}

Then we arrive to the following

\begin{proposition}\label{prop-8-3-10} 
Let in the described setup conditions \textup{(\ref{8-3-30})}, \textup{(\ref{8-3-31})} and \textup{(\ref{8-3-34})} be fulfilled.

Let $Q_1,Q_2$ be $h$-pseudo-differential operators with the symbols supported in $T^*X_2$ and 
\begin{equation}
\upmu_{2,\tau} (\Pi_{2,\tau}\cap \supp q_1\cap \supp q_2)=0
\label{8-3-35}
\end{equation}
where $\Pi_{2,\tau}$ is the set of point $(y,\eta)\in \Sigma_{2,\tau}$ periodic with respect to $\Psi_{2,t}$.

Then for $\Gamma (Q_{1x}e(.,.,\tau)\,^t\!Q_{1y})$ the standard Weyl two-term asymptotics holds with the remainder estimate $o(h^{1-d})$.
\end{proposition}

\begin{remark}\label{rem-8-3-11}
(i) Note that
\begin{equation}
\upmu_{j,\tau} (\Lambda_{j,\tau})=0
\label{8-3-36}
\end{equation}
where $\Lambda_{j,\tau}=\Lambda'_{j,\tau}\cup \Lambda''_{j,\tau}$,
$\Lambda'_{j,\tau}$ is the set of dead-end points of billiard (not generalized billiard) flow and $\Lambda''_{j,\tau}$ is the set of points $(y,\eta)\in \Sigma_{j,\tau}$ such that $\iota_j\Psi_{j,t}(y,\eta)\in \Xi_{3-j,\tau}$ where $\Xi_{k,\tau}=\iota_k \{(x,\xi)\in \Sigma_{k,\tau}, \{a_k,x_1\}=0$.

\medskip\noindent
(ii) Note that the grows properties of the branching Hamiltonian flow with reflections $\Psi_t$ coincide with those of $\Psi_{2,t}$ as long as we ensure that along billiards we do not approach to glancing rays on $X_1$:
\begin{multline}
\{a_1,x_1\}(x,\xi)\ge h^{\delta_1} \qquad 
\forall (x,\xi)\in \iota_1^{-1} \iota_2 \Psi_{2,t}(y,\eta)\cap \Sigma_{1,\tau}\\[3pt]
\forall (y,\eta)\in \supp q_2\ \forall t: \pm t\in [0,T]\  \Psi_{2,t}(y,\eta)\in T^*X_2|_Y
\label{8-3-37}
\end{multline}
where $T=T(h)$.
Then in assumptions of section~\ref{book_new-sect-7-4} one can recover similar results: namely that two-term Weyl asymptotics holds with the remainder estimate $O(T(h)^{-1}h^{1-d})$ which is usually $O(h^{1-d}|\log h|^{-1})$ or 
$O(h^{1-d+\delta})$ (depending on the growth and non-periodicity conditions to $\Psi_{2,t}$. Exact statement and proof are left to the reader (which is a relatively easy problem).

\medskip\noindent
(iii) In (ii) we need also to assume a more sharp version of (\ref{8-3-34}): namely if we consider approximate solutions to (\ref{8-3-32})--(\ref{8-3-33}) with precision $\varepsilon$ we get set $\Lambda_{\tau,\varepsilon}$ and we need to impose condition $\mes \Lambda_{\tau,\varepsilon}=o(\varepsilon^\delta)$ or $\mes \Lambda_{\tau,\varepsilon}=o(|\log\varepsilon|^{-1})$ as $\varepsilon\to +0$. 

\medskip\noindent
(iv) In subsubsection~``\nameref{book_new-sect-7-4-3-6}'' of subsection~\ref{book_new-sect-7-4-3} we did not consider points $z\in T^*Y$ where some of operators $a_j$ are elliptic on $\iota^{-1}z$ rather than hyperbolic. If we want to consider such points we need to assume (\ref{8-3-34}) or more sharp version of it described in (iii).

\medskip\noindent
(v) (\ref{8-3-37}) is guaranteed  if either
$a_2$ is cylindrically symmetric on $T^*X_2$ or if there no internal reflection for billiards from $X_2$ at all:
\begin{multline}
|\{a_1,x_1\} (x,\xi) |\ge \epsilon |\{a_2,x_1\}(x,\eta)|^l\\[3pt]
\forall (x,\xi)\in T^*X_1|_Y\cap\Sigma_{1,\tau}\
\forall (x,\eta)\in T^*X_2|_Y\cap\Sigma_{2,\tau}: \iota_1(x,\xi)=\iota_2(x,\eta).
\label{8-3-38}
\end{multline}
\end{remark}

\section{Analysis in $X_1$}
\label{sect-8-3-3-2}

However if operators $Q_1$ and $Q_2$ have symbols supported in $X_1$ situation is very different: trajectories of $a_1$ are periodic but the typical closed branching billiard  is the one which does not contain segments of Hamiltonian trajectories of $a_2$. More precisely, due to conditions
\begin{equation}
\upmu_{2,\tau} (\Pi_{2,\tau})=0
\label{8-3-39}
\end{equation}
and (\ref{8-3-36}) we conclude that 
\begin{equation}
\upmu_{1,\tau} (\Pi'_{1,\tau})=0
\label{8-3-40}
\end{equation}
where $\Pi'_{1,\tau}$ is the set of points  $z\in\Sigma_{1,\tau}$ such that there exists a billiard trajectory of $\Psi_t$ starting and ending in $z$ and which is not entirely billiard of $\Psi_{1,t}$.

Then for arbitrarily large $T>0$ and arbitrarily small $\varepsilon>0$ there exists a set $\cU=\cU_{T,\varepsilon}\subset T^*X$ such that 
\begin{gather}
\upmu_{1,\tau}(\Sigma_{1,\tau}\setminus \cU)\le \varepsilon,\label{8-3-41}\\[2pt]
\dist(\cU, \Pi'_{1,\tau,T}\cup \Lambda_{1,\tau,T})\ge \gamma\label{8-3-42}\\[2pt]
\dist(Y, \uppi_x \cU)\ge \gamma \label{8-3-43}
\end{gather}
with $\gamma=\gamma (T,\varepsilon)>0$ were $\Pi'_{1,\tau,T}$ s the set of points  $z\in\Sigma_{1,\tau}$ such that there exists a billiard trajectory of $\Psi_t$ of the length not exceeding $T$, starting and ending in $z$ and which is not entirely billiard of $\Psi_{1,t}$ and $\Lambda_{1,\tau,T}$ is the set of 
points  $z\in\Sigma_{1,\tau}$ such that billiard trajectory of $a_1$ hits the boundary at point of $\Lambda_{2,T,\tau}$. As usual, these sets increase as $T$ increases, their unions (with respect to $T$) coincide with $\Pi'_{1,\tau}$ and $\Lambda_{1,\tau}$ respectively and $\Lambda_{1,\tau,T}$ and $\Lambda_{1,\tau,T}\cap \Pi'_{1,\tau,T}$ are closed sets.

Therefore there exists operator 
$Q =Q_{T,\varepsilon}$ such that $\supp q\subset \cU$  and 
\begin{equation}
|\Gamma \bigl(Q_{1x}e(.,.,\tau)\,^t\!((I-Q)Q_2)_{y}\bigr)-
\kappa' _0 h^{-d}-\kappa' _1h^{1-d}|\le C_0\varepsilon h^{1-d}
\label{8-3-44}
\end{equation}
with $\kappa'_j=\kappa'_{j, Q_1,Q_2,Q}(\tau)$. 

This estimate holds for any fixed $\tau$ satisfying  (\ref{8-3-39}). Therefore we need to consider $\Gamma \bigl(Q_{1x}e(.,.,\tau)\,^t\!Q_{2y}\bigr)$ with $\supp q_2\subset \cU$.

We can assume without any loss of the generality that  (\ref{8-3-42}), (\ref{8-3-43}) are fulfilled for all $\tau$. Then for $\epsilon \le | t| \le T$  
\begin{equation}
\Gamma \bigl(Q_{1x}U(.,.,t)\,^t\!Q_{2y}\bigr) \equiv 
\Gamma \bigl(Q_{1x}Z(.,.,t)\,^t\!Q_{2y}\bigr)
\label{8-3-45}
\end{equation}
where as $\pm t>0$\ \ $Z(t)$ is Schwartz kernel of ``approximate propagation with scattering semigroup'' $\mathbf{Z}(t)$ constructed in the following way (for a sake of simplicity in notations we consider $t>0$: 

\medskip\noindent
(i) We represent $Q_2=Q_{2,1}+\dots + Q_{2,N}$ where symbols $Q_{2,\nu}$ have small supports;

\medskip\noindent
(ii) For each $\nu$ we select $0=t_{0,\nu} < \dots < t_{M,\nu}=t$ in such a way that for any $z\in \supp q_{2,\nu}$ \ $\Psi_{t'} (z)$ is disjoint from $Y$ as $t'=t_{j,\nu}$ and  $\Psi_{t'} (z)$ hits $Y$ no more than once as 
$t_{j,\nu}\le t'\le t_{j+1,\nu}$;

\medskip\noindent
(iii) We set
\begin{equation*}
\mathbf{Z}(t)Q_{2,\nu}= \Bigl(\prod_{0\le j\le M-1} 
\psi e^{ih^{-1} (t_{j+1,\nu}- t_{j,\nu})f(A)}\Bigr) \cdot Q_{2,\nu}
\end{equation*}
where $\psi \in \sC^\infty$ is supported in 
$\{x\in X_1, \dist (x,Y)\ge \gamma'\}$ and equals $1$ in 
$\{x\in X_1, \dist (x,Y)\ge 2\gamma'\}$.

Then 
\begin{gather}
\mathbf{Z}(t_1+t_2) Q_2\equiv \mathbf{Z}(t_2) \mathbf{Z}(t_1) Q_2
\qquad\text{as\ \ } \pm t_1>0, \pm t_2>0\label{8-3-46}\\
\shortintertext{and}
\mathbf{Z}(t)^*\equiv \mathbf{Z}(-t)\label{8-3-47}
\end{gather}
and due to periodicity of $\Psi_{1,t}$ with period $T_0$  $\mathbf{Z}(\pm T_0)$ are $h$-pseudo-differential operators.

\begin{remark}\label{rem-8-3-12}
In contrast to what we had before they are not necessarily unitary because \emph{reflection coefficients\/} are not unitary anymore. More precisely, 

\medskip\noindent
(i) If $\Psi_{1,t}$ hits $T^*X|_Y$ at points $(x',\xi)$ with $(x',\xi')\notin \iota_2 \Sigma_{2,\tau}$ reflection coefficient $\varkappa_{11}(x',\xi',\tau)$ still satisfies $|\varkappa_{11}(x',\xi',\tau)|=1$.

\medskip\noindent
(ii) However if
$\Psi_{1,t}$ hits $T^*X|_Y$ at points $(x',\xi)$ with 
$(x',\xi')\in \iota_2\Sigma_{2,\tau}$ then reflection-refraction matrix
$\varkappa(x',\xi',\tau)=\bigl(\varkappa_{jk}(x',\xi',\tau)\bigr)_{j,k=1,2}$ is unitary but reflection coefficient $\varkappa_{11}(x',\xi',\tau)$ satisfies only $|\varkappa_{11}(x',\xi',\tau)|\le 1$ and the equality means  exactly that $\varkappa_{12}(x',\xi',\tau)=\varkappa_{21}(x',\xi',\tau)=0$ (then $|\varkappa_{22}(x',\xi',\tau)|=1$ as well).
\end{remark}

So the principal symbol of $\mathbf{Z}(T_0)$ is the product of $e^{i\ell_0(x,\xi)}$ where $\ell_0(x,\xi)$ is given by the usual formula for manifolds without boundary (but with trajectories replaced by billiards) and also of reflection coefficients coefficients in the point of reflection. If the reflection coefficients do not vanish we can rewrite the answer in the form
\begin{equation}
\mathbf{Z}(T_0)\equiv e^{i\varepsilon h^{-1}L},\qquad 
\mathbf{Z}(-\varepsilon h^{-1} T_0)\equiv e^{-iL^*}
\label{8-3-48}
\end{equation}
with $h$-pseudo-differential operator $L$, $L+L^* \le 0$ and $\varepsilon =h$ (but in calculations below we are not assuming this). 

We will use this formula even if the coefficients vanish: in this case we simply formally allow $\Im \ell = +\infty$.

Then as $\bar{\chi}\in \sC^\infty_0 [\epsilon -1, 1-\epsilon ]$
\begin{multline}
F_{t\to h^{-1}\tau}
\bar{\chi}_{T_0}  (t-nT_0) \Gamma \bigl(Q_{1x}Z(.,.,t)\,^t\!Q_{2y}\bigr)= \\h^{1-d}J(h,\varepsilon, n,\tau)\cdot \widehat{\bar{\chi}}\bigl(\frac{\tau}{h}\bigr)
+O\bigl(h^{2-d}\bigr),\label{8-3-49}
\end{multline}
with
\begin{equation}
J(h,\varepsilon, n,\tau) = (2\pi )^{1-d}\int _{\Sigma_{1,\tau}}  e^{i h^{-1}((\varepsilon \Re \ell -\tau )n +\varepsilon \Im \ell |n|)}\, d\upmu_{1,\tau}.
\label{8-3-50}
\end{equation}
We can see easily that
\begin{claim}\label{8-3-51}
$J(h,\varepsilon, n,\tau)=o(1)$ as $n\to \infty$ provided
\begin{equation}
\upmu_\tau \bigl(\bigl\{(x,\xi)\in \Sigma_\tau:\ 
\Im \ell (x,\xi)= \nabla \Re \ell (x,\xi)=0\bigr\}\bigr)=0.
\label{8-3-52}
\end{equation}
\end{claim}
On the other hand, note that
\begin{equation}
 \sum_{n\in \bZ\setminus 0}
e^{i h^{-1}(n\alpha  -|n|\beta)}=
\frac{2 e^{-\beta}\cos\alpha  -2 e^{-2\beta}}
{1- 2 e^{-\beta}\cos \alpha + e^{-2\beta}}\qquad\text{as\ \ } \beta>0\label{8-3-53}
\end{equation}
and formally we can extend this for $\beta=0$ as well. 

Then after summation with respect to $n$ for $\beta=0$ we get exactly $\partial_\alpha \Upsilon  (\alpha)$ with $\Upsilon$ defined by (\ref{book_new-6-2-63}) and for $\beta>0$ we get that the expression (\ref{8-3-53}) is equal to
$\partial_\beta  \Upsilon *_\alpha (2\pi)^{-1} G(\alpha,\beta)$ where 
$G(\alpha,\beta)=2\beta\bigl(\alpha^2+\beta^2\bigr)^{-1}$.

Note that 
\begin{claim}\label{8-3-54}
Function
\begin{equation}
\Upsilon (\alpha ,\beta)\Def  \Upsilon *_\alpha  G(\alpha,\beta)
\label{8-3-55}
\end{equation}
is harmonic function  as $\beta>0$, coincides with $\Upsilon(\alpha)$ as $\beta=0$, is $o(1)$ as $\beta\to +\infty$ and is periodic with respect to $\alpha$.
\end{claim}

Taking in account results of this subsection we arrive to 

\begin{theorem}\label{thm-8-3-13}
Let in the described setup conditions \textup{(\ref{8-3-30})}, \textup{(\ref{8-3-31})}, \textup{(\ref{8-3-34})}, \textup{(\ref{8-3-39})} and \textup{(\ref{8-3-52})} be fulfilled.

Then asymptotics
\begin{multline}
\Gamma e(.,.,\tau) 
=\kappa_0(\tau) h^{-d} +\\
\Bigl( \kappa_1(\tau) + \int _{\Sigma_\tau}
\Upsilon \bigl(h^{-1}(\varepsilon\Re\ell-\tau),-h^{-1}\varepsilon 
\Im \ell\bigr)\,d\upmu_\tau \bigr)\,d\upmu_\tau\Bigr)+o(h^{1-d})
\label{8-3-56}
\end{multline}
holds.
\end{theorem}

\begin{problem}\label{prob-8-3-14}
Calculate all coefficients $\varkappa_{jk}$ for two Schr\"odinger operators $A_1$ and $A_2$.
\end{problem}

We will solve this problem in the very special case:

\begin{example}\label{ex-8-3-15}
Assume that (at the given point $z\in T^*Y$) $a_j=c_j^2|\xi|^2$ and the boundary condition is
\begin{gather}
u_2=\alpha  u_1,\label{8-3-57}\\
D_1u_2=-\beta  D_1u_1\label{8-3-58}
\end{gather}
(recall that $x_1>0$ in both $X_1,X_2$; in more standard notations of $x_1<0$ in $X_2$ one needs to skip ``$-$''). Then $\{A,B\}$ to be self-adjoint requires $\alpha\beta^\dag= c_1^2c_2^{-2}$. We do not consider gliding points since their measure is $0$.

\medskip\noindent
(i) Consider point $\xi'$ with $c_1^2 |\xi'|^2< \tau <c_2^2|\xi'|^2$. One can prove easily that
\begin{gather}
\varkappa_{11}(\xi',\tau)=e^{2i\varphi}\label{8-3-59}\\
\shortintertext{with}
\varphi (\xi',\tau)= \arctan \Bigl( |\beta|^{-2}c_2^{-2}c_1^2\bigl(c_2^2|\xi'|^2-\tau\bigr)^{\frac{1}{2}}
\bigl(\tau -c_1^2|\xi'|\bigr)^{-\frac{1}{2}}\Bigr)\label{8-3-60}
\end{gather}
and therefore assumption (\ref{8-3-52}) is fulfilled provided $c_1\ne c_2$, $c_j,\alpha,\beta$ are symmetric (i.e. $c_j\varrho=c_j$) and the subprincipal symbol is $0$. If $c_1>c_2$ this case does not appear.

\medskip\noindent
(ii) Consider point $\xi'$ with $c_j^2 |\xi'|^2< \tau$ for both $j=1,2$.  One can prove easily that 
\begin{gather}
\varkappa_{11}(\xi',\tau)=\bigl(\omega -1\bigr)\bigl(\omega +1\bigr)^{-1} \label{8-3-61}\\
\shortintertext{with}
\omega (\xi',\tau)= |\beta|^{-2}c_2^{-2}c_1^2\bigl(\tau-c_2^2|\xi'|^2\bigr)^{\frac{1}{2}}
\bigl(\tau -c_1^2|\xi'|\bigr)^{-\frac{1}{2}}\label{8-3-62}
\end{gather}
and $|\varkappa_{11}|<1$ so assumption (\ref{8-3-52}) is fulfilled.
\end{example}

The following problems seem to be rather straightforward but not easy and worth publishing:
\begin{problem}\label{prob-8-3-16}
Under stronger conditions to $\Psi_{2,t}$ and $\ell$ (see (\ref{8-3-20})) prove remainder estimate $O(h^{1-d+\delta})$.
\end{problem}

\begin{problem}\label{prob-8-3-17}
Prove the same results as Hamiltonian billiards of $a_1$ are assumed to be periodic only on one energy level $\tau$ rather than on neighboring levels.
\end{problem}

\chapter{Two periodic flows}
\label{sect-8-3-4}

Now we consider the case of both flows $\Psi_{1t}$ and $\Psi_{2t}$ generated by $f(a_1)$ and $f(a_2)$ being periodic and satisfying (\ref{8-3-30}) as in the following figures:
\begin{figure}[h!]\centering
\subfloat[Hemisphere, projected to equatorial plane]{\includegraphics[width=.5\textwidth]{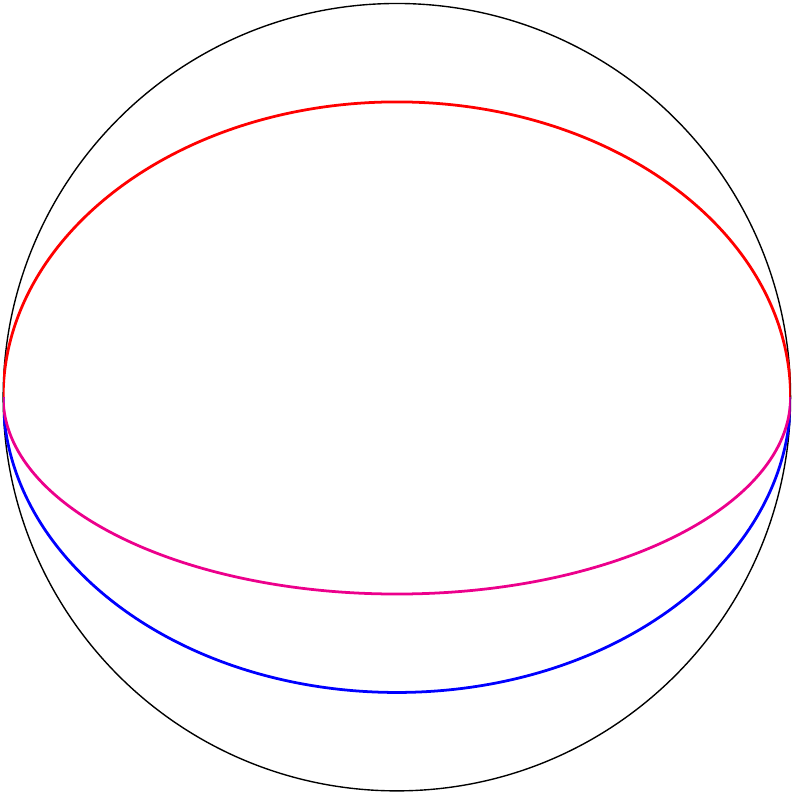}}
\qquad
\subfloat[half-plane]{\includegraphics[width=.4\textwidth]{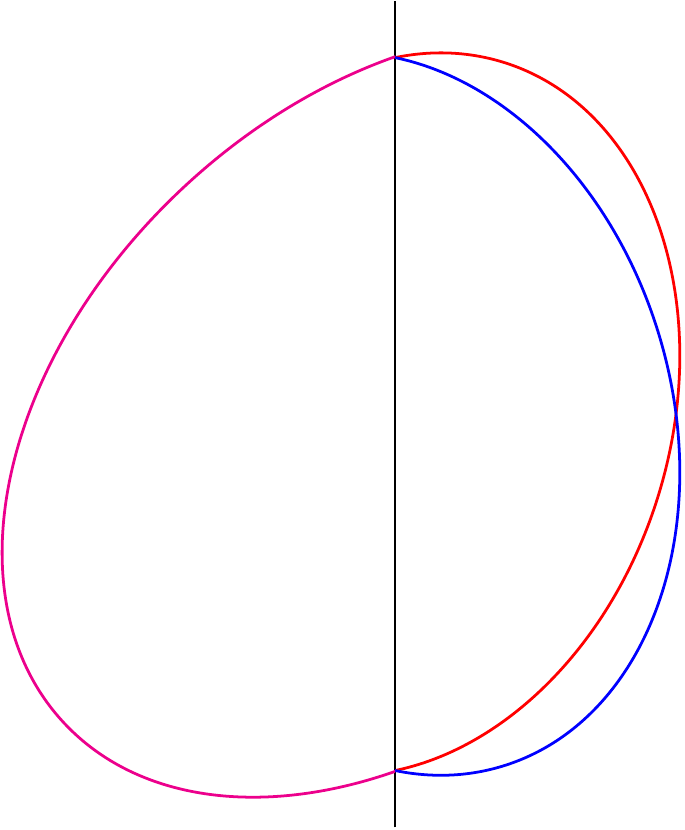}} \caption{\label{fig-branch-billiards} $x$-projection of billiard trajectory on two half-planes, glued together: original (red), reflected (blue), refracted (magenta)}
\end{figure}

\section{Examples and discussion}
\label{sect-8-3-4-1}

Let us start from an example:

\begin{example}\label{ex-8-3-18}
Consider two Schr\"odinger operators in $\bR^d_+= \bR^+ \times \bR^{d-1}$ with the principal symbols
\begin{equation}
a_j (x,\xi)= \frac{1}{2}\omega_j^2 \bigl(|\xi|^2+|x|^2|- E_j\bigr)
\label{8-3-63}
\end{equation}
intertwined through boundary conditions. Then $T_j= 2\pi \omega_j^{-2}$.

We are interested in the number $\N^-_h(0)$ of the negative eigenvalues of operator $A_B$. In order to consider the eigenvalue counting function $\N^-_h(\lambda)$ for arbitrary  spectral parameter $\lambda$ one needs to redefine
\begin{equation}
E_j\Def E_j+2\omega_j^2 \lambda.
\label{8-3-64}
\end{equation}

But $\N^-(0)$ is also the maximal dimension of the \emph{negative subspace\/}
of the operator $A$ which does not change if we replace $A$ by $J^{-\frac{1}{2}}AJ^{-\frac{1}{2}}$ with the positive  self-adjoint operator $J$ and the latter equals to $\N^-_h(0)$ for the latter operator. Picking up $J=J(x)$ equal $\omega_j^2$ in $X_j$ we arrive to the case of $\omega_j=1$ and $T_j=2\pi$. Therefore in this example one can make periods equal $T_1=T_2$. The dependence of $\omega_j$ will come through $E_j$ redefined by (\ref{8-3-64}).
\end{example}

There are however many deficiencies of the described approach; most important is that it does not work with the following example (unless some unnatural conditions to $\omega_j$, $E_j$ are imposed):

\begin{example}\label{ex-8-3-19}
Consider two  operators on hemisphere $\bS^d_+$ with the principal symbols
\begin{equation}
a_j (x,\xi)= \omega_j^2\bigl(|\xi|^2 -E_j\bigr)
\label{8-3-65}
\end{equation}
where $|\xi|$ is calculated in the standard Riemannian metrics.
\end{example}

The main obstacle here is that multiplying operator by $J^{-\frac{1}{2}}$ from both sides and calculating $f(A)$ do not play well together.

\section{Reduction to the boundary}
\label{sect-8-3-4-2}

First of all notice that the contribution to the remainder of the domain 
$\{z, \dist (z,T^*Y)\le \varepsilon_0\}$ does not exceed 
$C_0\epsilon_0^2 h^{1-d}$ and therefore one needs to consider the contribution of $\{z, \dist (z,T^*Y)\ge \varepsilon_0\}$ with arbitrarily small but fixed $\varepsilon_0$.

In this zone one can consider $U(x,y,t)$ as the solution of the Cauchy problem with respect to $x_1$ (or $y_1$) with data at $\{x_1=0\}$ (or $\{y_1=0\}$). Rewriting therefore $U = R_{x}U\,^t\!R_{y}$ where $R$ is an operator resolving this Cauchy problem is given by an oscillatory integral we can rewrite 
\begin{equation}
\Gamma \bigl(U(.,.,t)\,^t\!Q_{y}\bigr) \equiv
\Gamma ' \bigl( \cQ(x',hD'_{x},hD_{t},h) V\bigr)
\label{8-3-66}
\end{equation}
with $h$-pseudo-differential operator $R$ and $V=\eth_{x,m-1}\eth_{y,m-1}U$ the Cauchy data for $U$ (assuming that $m$ is an order of $A$). 

One can see easily that as $Q$ is an operator with the symbol supported in $T^*X_j$, the principal symbol of $\cQ$ at $(x',\xi,\tau)$ is defined as an averaging of $Q$ along $\Psi_{jt}(z_j)$ with $t\in [0, t_j(z_j)]$ where 
$z_j= \iota^{-1}(x',\xi') \cap \Sigma_{j\tau}$ and $t_j(z_j)$ is the time of the next hit of the boundary.

Then instead of continuous family of Fourier integral operators $e^{ih^{-1}tA}$ we can consider a discrete family $\cF^n$ defined by the following way: consider 
$v_j=\eth_{m-1}u_j = (v_j^-,v_j^+)$ with $v_j^\mp$ corresponding to incoming and outgoing solutions. 

We have two types of trajectories: those in $\Sigma_{j\tau}$ which hit $T^*X|_Y$ at points elliptic for $(a_{3-j}-\tau)$ where $j=1,2$ and those which hit $T^*X|_Y$ at points hyperbolic for $(a_k-\tau)$ for both $k=1,2$. The analysis of the former is of no different from what we have seen in the previous subsection: we just consider $\Psi_{jt}$ and construct the (real-valued) symbol $\ell$ as long as we assume that $(x',\xi)\notin \Lambda$ where $\Lambda$ is defined in the paragraph preceding (\ref{8-3-34}).

The analysis of the latter is more interesting. Therefore as before we have a  matrix of reflection-refraction $(\varkappa_{jk})_{j,k=1,2}$ which is a symbol defined on $\cU\times (\tau-\varepsilon',\tau+\varepsilon')$ where $\cU$ is a a zone described above. This matrix is unitary in the norm
\begin{equation}
\|V\|= \Bigl(\sum_j \bigl(|V^+_j|^2 + |V^-_j|^2\bigr)|\{a_j,x_1\}|_{x_1=0}\Bigr)^{\frac{1}{2}}.
\label{8-3-67}
\end{equation}

We also have diagonal unitary matrices  
$M'_\nu =\diag (e^{i\ell_{1\nu}}, e^{i\ell_{2\nu}})$ with $\ell_{j\nu}$ calculated along $\nu$-th leg of the closed trajectory of $\Psi_{jt}$ with $t=\frac{1}{2}T_0$ and $\nu=1,2$.

Consider closed branching billiard originated from point $z\in T^*Y$; first there are two trajectories (in $X_1$ and $X_2$) both hitting $T^*X|_Y$ at $\varrho(z)$ and then there are two trajectories (in $X_1$ and $X_2$)  returning to $z$.  Let us recall that $\varrho(z)$ is an antipodal point.

Therefore branching billiard is characterized by an unitary (in (\ref{8-3-67})-norm) matrix 
\begin{equation}
M(z,\tau)= \varkappa (z ,\tau) \,M' (z,\tau)\, \varkappa (\varrho(z) ,\tau) \,M' (\varrho(z),\tau).
\label{8-3-68}
\end{equation}

\section{Analysis of the evolution}
\label{sect-8-3-4-3}
One can prove easily that modulo $O_T(h^{2-d})$
\begin{multline}
F_{t\to h^{-1}\tau} \chi_T (t)
\Gamma ' \bigl( \cQ(x',hD'_{x},hD_{t},h) V\bigr)\equiv\\[3pt]
(2\pi h)^{1-d} \chi_T (T_1+\dots+ T_n)
\iint \sum_{N\ne 1}\sum_{j=(j_1,\dots,j_N)\in \{1,2\}^N}\\[3pt]
M_{j_1j_2}M_{j_2j_3}\cdots M_{j_{N-1}j_N} e^{-ih^{-1}\tau (T_{j_1}+\dots T_{j_N})} \cQ_{j_Nj_1}  dx' d\xi'
\label{8-3-69}
\end{multline}
where $\chi\in \sC_0^\infty $ is supported in $[\frac{1}{2},1]$ and symbols $M_{jk}$ and $\cQ$ and $T_{j}$ are calculated at $(x',\xi',\tau)$ where $T_j=T_j(\tau)$ since all trajectories are periodic on energy levels close to $\tau$.

One can rewrite terms with equal $N$ as
\begin{multline}
(2\pi h)^{1-d} \int\hat{\chi}(\lambda)e^{iT^{-1} \lambda (T_{j_1}+\dots+ T_{j_N})}
\iint \sum_{j=(j_1,\dots,j_N)\in \{1,2\}^N}\\[3pt]
M_{j_1j_2}M_{j_2j_3}\cdots M_{j_{N-1}j_N} e^{-ih^{-1}\tau (T_{j_1}+\dots T_{j_N})} \cQ_{j_Nj_1}  dx' d\xi'd\lambda
\label{8-3-70}
\end{multline}
which in turn equals
\begin{equation}
(2\pi h)^{1-d} \int\hat{\chi}(\lambda)
\iint\\
\tr \bigl(S(\lambda,\tau,x',\xi')^N \cQ\bigr)  dx' d\xi'd\lambda
\label{8-3-71}
\end{equation}
where 
\begin{equation}
S(\lambda T^{-1},\tau,x',\xi')=
\begin{pmatrix} e^{i(\lambda T^{-1}+ \tau h^{-1})T_1} m_{11} &
e^{i(\lambda T^{-1}+ \tau h^{-1})T_1} m_{12}\\
e^{i(\lambda T^{-1}+ \tau h^{-1})T_2} m_{21} &
e^{i(\lambda T^{-1}+ \tau h^{-1})T_2} m_{22}
\label{8-3-72}
\end{pmatrix}
\end{equation}
and $M=(m_{jk})$.

Consider eigenvalues of this matrix; as $T\to \infty$ they tend to eigenvalues of the same matrix with $\mu=0$. If 
\begin{multline}
\mes \{(x',\xi'): \rho \in \Spec \bigl(S(0,\tau,x',\xi')\bigr)\}=o(1)\\[2pt]
\forall \tau \in 
[\bar{\tau}-\epsilon hT^{-1}, \bar{\tau}+\epsilon hT^{-1}] \qquad
\text{as \ \ }T\to \infty 
\label{8-3-73}
\end{multline}
then this term is $o(h^{1-d})$ and therefore we estimated expression (\ref{8-3-69}) by $o(Th^{1-d})$ which in the end of the day returns remainder estimate $o(h^{1-d})$ with the Tauberian main part. 

Note that instead of $S$ we can consider matrix 
\begin{equation}
\tilde{S}(\tau,x',\xi')=
\begin{pmatrix} e^{i\tau h^{-1} T^*} m_{11} &
 m_{12}\\[2pt]
 m_{21} &
e^{-i\tau h^{-1} T^*} g_{22}
\label{8-3-74}
\end{pmatrix}
\end{equation}
with $T^*=\frac{1}{2}(T_1-T_2)$ which is different from $S(0, .,.,.)$ by factor
$e^{\frac{1}{2}ih^{-1}\tau (T_1+T_2)}$.

Consider correction to the main part; we are interested at $\tau=0$; so we integrate (\ref{8-3-69}) with $\chi=1$ from $\tau=-\infty$ to $\tau=0$ resulting in (after multiplication by $h^{-1}$ and modulo $o(h^{1-d})$)
\begin{multline}
(2\pi h)^{1-d} \chi_T (T_1+\dots+ T_n)
\iint \sum_{N\ne  1}\sum_{j\in \{1,2\}^N}\\
M_{j_1j_2}M_{j_2j_3}\cdots M_{j_{N-1}j_N}  (T_{j_1}+\dots T_{j_N})^{-1} \cQ_{j_Nj_1}  dx' d\xi'
\label{8-3-75}
\end{multline}
where we plug $\tau=0$.

This latter equals to
\begin{equation}
\Omega (\tau)= (2\pi h)^{1-d}\iint
\tr \Upsilon (L(x',\xi', h^{-1}\tau))\,\cQ(x',\xi.\tau) dx'd\xi'
\label{8-3-76}
\end{equation}
where $L$ is defined by $e^{iL}=S$ and $\Upsilon$ is defined by (\ref{book_new-6-2-63}).

\begin{theorem}\label{thm-8-3-20}
Consider two Schr\"odinger operators $A_1$ and $A_2$ with all periodic trajectories on levels close to $\tau$, satisfying \textup{(\ref{8-3-30})}.
Further, let \textup{(\ref{8-3-34})} and \textup{(\ref{8-3-73})} be fulfilled. Then asymptotics
\begin{equation}
\tr \Gamma \bigl(e(.,.,\tau)\bigr)=
\varkappa_0(\tau)h^{-d}+\bigl(\varkappa_1(\tau)+\Omega(\tau)\bigr)h^{1-d}+
o(h^{1-d})
\label{8-3-77}
\end{equation}
holds.
\end{theorem}

\begin{example}\label{ex-8-3-21}
Consider example~\ref{ex-8-3-15} in the new conditions to Hamiltonian flows. We can use results of that example immediately to treat billiards with complete internal reflections (thus non-branching).

To treat branching billiards we need to calculate eigenvalues of the matrix $(\varkappa_{jk})$. It follows from calculations of   example~\ref{ex-8-3-15} that $\varkappa_{22}=-\varkappa_{11}$; thus  we arrive to
\begin{equation}
\lambda = e^{\pm i \varphi},\qquad \varphi= \arccos \varkappa_{11}=
\arccos \bigl( (\omega-1)(\omega+1)^{-1};
\label{8-3-78}
\end{equation}
then obviously (\ref{8-3-73}) is fulfilled and asymptotics (\ref{8-3-77}) holds.
\end{example}

\begin{remark}\label{rem-8-3-22}
On the contrary assume that (\ref{8-3-73}) is not fulfilled. Then there will be eigenvalues of high multiplicity or clusters of eigenvalues located in $o(h)$-vicinities of solutions $\tau$ of equation $\det (S-1)=0$. One can rewrite this equation as 
\begin{equation}
\cos \bigl( -\varphi + \frac{1}{2}\tau h^{-1} (T_1+T_2)\bigr)=
\cos \alpha \cdot \cos \bigl( -\psi + \frac{1}{2}\tau h^{-1} (T_1-T_2)\bigr)
\label{8-3-79}
\end{equation}
as $M=\begin{pmatrix}\ \ e^{i(\phi+\psi)}\cos\alpha & e^{i(\phi+\chi)}\sin\alpha\\
-e^{i(\phi-\chi)}\sin\alpha & e^{i(\phi-\psi)}\cos\alpha\end{pmatrix}$
which is the general form of the unitary matrix.
\end{remark}

The following problem seems to be a difficult one and worth of publication:

\begin{problem}\label{prob-8-3-23}
Prove the same results as flows $\Psi_{j,t}$ are assumed to be periodic only on one energy level $\tau$ rather than on neighboring levels. 
\end{problem}

On the contrary, the following problem seems to be relatively easy:

\begin{problem}\label{prob-8-3-24}
(i) Prove the same results as there are more than two manifolds $X_j$ $j=1,\dots,m'$ with the all billiards periodic.

\medskip\noindent
(ii)  Extend these results to the case when there are also manifolds $X_j$ $j=m'+1,\dots,m$ with almost all billiards non-periodic.
\end{problem}

\input IRO24.bbl

\end{document}

%% file: defines_IRO24.tex
\theoremstyle{plain}
\newtheorem{theorem}{Theorem}[chapter]

\newtheorem{proposition}[theorem]{Proposition}
\newtheorem{corollary}[theorem]{Corollary}

\theoremstyle{definition}

\theoremstyle{remark}
\newtheorem{remark}[theorem]{Remark}
\newtheorem{example}[theorem]{Example}
\newtheorem{problem}[theorem]{Problem}

\numberwithin{equation}{chapter}

\DeclareMathAlphabet{\mathpzc}{OT1}{pzc}{m}{it}

 \newcommand{\cF}{\mathcal{F}}

 \newcommand{\cQ}{\mathcal{Q}}

 \newcommand{\cU}{\mathcal{U}}

 \newcommand{\sC}{\mathscr{C}}

 \newcommand{\sH}{\mathscr{H}}

\newcommand{\N}{{\mathsf{N}}}

\newcommand{\const}{{\mathsf{const}}}
\newcommand{\dist}{{{\mathsf{dist}}}}

\newcommand{\bR}{{\mathbb{R}}}
\newcommand{\bS}{{\mathbb{S}}}

\newcommand{\bZ}{{\mathbb{Z}}}

\newcommand{\fD}{{\mathfrak{D}}}

\def\1{\boldsymbol {|}}
\newcommand{\Def}{\mathrel{\mathop:}=}

\newcommand{\diag}{\operatorname{diag}}

\renewcommand{\Im}{\operatorname{Im}}       % \Im is already frak I

\newcommand{\mes}{\operatorname{mes}}

\renewcommand{\Re}{\operatorname{Re}}       % \Re is already frak R

\newcommand{\Spec}{\operatorname{Spec}}

\newcommand{\supp}{\operatorname{supp}}

\newcommand{\tr}{\operatorname{tr}}
\newcommand{\Tr}{\operatorname{Tr}}

\newenvironment{claim}[1][{\textup{(\theequation)}}]{\refstepcounter{equation}\vglue10pt
\begin{trivlist}
\item[{\hskip\labelsep#1}]}{\vglue10pt\end{trivlist}}

\newenvironment{claim*}[1][{}]{\vglue10pt
\begin{trivlist}
\item[{\hskip\labelsep#1}]}{\vglue10pt\end{trivlist}}

\newcounter{note}

\DeclareTextCommand{\textinfty}{PU}{\9042\036}

\DeclareTextCommand{\textge}{PU}{\9042\145}
\DeclareTextCommand{\textle}{PU}{\9042\144}
\DeclareTextCommand{\texthat}{PD1}{\136}

%% file: IRO24.bbl
\bibliographystyle{alpha}

\providecommand{\bysame}{\leavevmode\hbox to3em{\hrulefill}\thinspace}

\vglue .06truein

\begin{tabular}{rrl}
&{\hskip 200 pt} &Department of Mathematics,\cr
&&University of Toronto,\cr
&&40, St.George Str.,\cr
&&Toronto, Ontario M5S 2E4\cr
&&Canada\cr
&&ivrii@math.toronto.edu\cr
&&Fax: (416)978-4107\cr
\end{tabular}